\documentclass[11pt]{amsart}
\title{Hamiltonian Actions on Homogeneous Bounded Domains}
\author{Maxim Kukol}
\address{School of Mathematics and Natural Sciences, University of Wuppertal, Gaußstr. 20, 42119 Wuppertal, Germany}
\email{kukol@uni-wuppertal.de}
\thanks{This work was partially supported by the SFB TRR 191 {\it Symplectic Structures in Geometry, Algebra and Dynamics}, funded by the DFG (Projektnummer 281071066~--~TRR 191).}

\usepackage[utf8]{inputenc}
\usepackage[english]{babel}
\usepackage{amsmath, amsfonts, amssymb, mathtools, mathdots, faktor, polynom}

\usepackage{tabularx, caption, setspace, calc, multirow, bigdelim}

\usepackage{graphicx, pdfpages}
\usepackage{tikz}
\usepackage{tikz-cd}
\usetikzlibrary{arrows, fadings, patterns}

\usepackage{pgfplots}
\pgfplotsset{compat=1.18}

\usepackage{hyperref}
\usepackage[english]{cleveref}

\usepackage{extarrows}
\usepackage{enumerate}

\theoremstyle{plain} % für Definitionen, Beispiele
\newtheorem{Def}{Definition}[subsection]
\newtheorem{Lem}[Def]{Lemma}
\newtheorem{Coro}[Def]{Corollary}
\newtheorem{Prop}[Def]{Proposition}
\newtheorem{Theo}[Def]{Theorem}

\theoremstyle{definition} % für aufrechte Schrift
\newtheorem{Ex}[Def]{Example}

\newtheorem{Rk}[Def]{Remark}

\theoremstyle{plain} % für Theoreme
\newtheorem{THEO}{Theorem}

\newtheorem*{Deff}{Definition}

\newcommand{\Proof}{\textit{Proof.}\ }

\newcommand{\NN}{\mathbb N}

\newcommand{\RR}{\mathbb R}
\newcommand{\CC}{\mathbb C}

\newcommand{\MM}{\mathcal M}

\DeclareMathOperator{\im}{Im}

\begin{document}
\maketitle
\begin{abstract}\vspace{-3em}\emergencystretch=1em
	We investigate Hamiltonian actions of non-compact Lie groups on a homogeneous bounded domain $X$. As a main result, we point out a Lie-theoretical condition for a closed Lie group $H$ of the automorphism group of $X$ which ensures that the symplectic reduction $\mu^{-1}(0)/H$ with respect to the momentum map $\mu$ at hand, is a Stein manifold. Moreover, for the class of connected subgroups of translations the quotient $(H^\mathbb{C}\cdot X)/H^\mathbb{C}$ is realized as a Siegel domain and we show that the symplectic reduction $\mu^{-1}(0)/H$ is biholomorphic to such a Stein quotient.
\end{abstract}

\section{Introduction}
Let $X$ be a Kähler manifold with a compact Lie group $K$ acting on $X$ by holomorphic transformations in a Hamiltonian fashion. If $X$ is Stein and the complexification $K^\CC$ acts on $X$ holomorphically, there is a well-developed geometric invariant theory showing that if the momentum map is induced by a $K$-invariant potential, the symplectic reduction by $K$ coincides with the categorical quotient $X//K^\CC$. Thus the symplectic reduction naturally carries the structure of a Stein space (see \cite{GS}, \cite{HHL} or \cite{HL}). For example, this is an essential ingredient in proving an equivariant version of Grauert's Oka Principle (see \cite{HK}). Moreover, for real Lie subgroups of $K^\CC$ which are compatible with the Cartan decomposition of $K^\CC$, analogous results were obtained in \cite{HG} using a modified momentum map. However, if the complexification does not act, the situation is much more complicated (see \cite{H}).\vspace{6pt}
	
Let $X$ be a bounded domain in $\CC^d$ and $\textup{Aut}_{\mathcal{O}}(X)$ the group of holomorphic automorphisms of $X$. The group $\textup{Aut}_{\mathcal{O}}(X)$ is a real Lie group which acts properly on $X$ (see \cite{C}). Assume that $X$ is homogeneous and let $G$ be the connected component of $\textup{Aut}_{\mathcal{O}}(X)$ which contains the identity. Since $X$ is homogeneous and connected, $G$ acts transitively on $X$. For a fixed base point $x_0\in X$, the isotropy group $G_{x_0}$ is a maximal compact subgroup of $G$ and there exists a decomposition $G=BG_{x_0}$, where $B$ is a simply connected split solvable subgroup of $G$ (see \cite{Kan} and \cite{Vin}). In particular, $B$ acts freely and transitively on $X$ (see \cite{Kan}, p. 4). Further, it is known that every homogeneous bounded domain is biholomorphic to a Siegel domain of the first or second kind (see \cite{P}, Appendix). On this Siegel domain the group $B$ acts by affine transformations (see \cite{Kan}, Ch. 1) and the action extends to the whole $\CC^d$. Moreover, the complexification $B^\CC$ acts on $\CC^d$ holomorphically.\vspace{6pt}
	
The Bergman metric on $X$ defines a Kähler form $\omega_\mathcal{B}$ on $X$. The form $\omega_\mathcal{B}$ is invariant under biholomorphic maps of $X$. In fact, the group $G$ acts in a Hamiltonian fashion on $X$. More precisely, there exists a $G$-equivariant real analytic map $\mu:X\longrightarrow\mathfrak{g}^*$, where $\mathfrak{g}^*$ denotes the dual vector space of the Lie algebra of $G$, such that for $x\in X$ and $\xi\in\mathfrak{g}$ 
    \begin{equation}\tag{$*$}\label{eq:mom}
	\textup{d}\mu^\xi(x)=\omega_\mathcal{B}(\xi_{X}(x),\cdot),
	\end{equation}
where $\xi_{X}(x):=\frac{d}{dt}\big|_{t=0}\exp(t\xi)\cdot x$ and $\mu^\xi(x):=\mu(x)(\xi)$. The map $\mu$ is called a momentum map. From the symplectic point of view, $\mu^\xi$ is a Hamiltonian function with Hamiltonian vector field $\xi_X$. For every closed subgroup $H$ of $G$ we have a restricted momentum map $\mu_H:X\longrightarrow\mathfrak{h}^*$. We denote by $\MM_H:=\mu_H^{-1}(0)$ the zero level set of $\mu_H$. By the equivariance of $\mu_H$, the group $H$ acts on $\MM_H$ and the quotient $\MM_H/H$ is the corresponding symplectic reduction of $X$.\vspace{6pt}
	
In this paper we assume that $H$ is a closed subgroup of $B$. It follows that $H$ acts freely and properly on $X$ and (\ref{eq:mom}) implies that $\MM_H$ is a smooth manifold. Moreover, the symplectic reduction of $X$ is a symplectic manifold whose symplectic structure $\widetilde\omega$ is characterised by 
    $$
	\pi_{\MM_H}^*\widetilde\omega=\iota^*\omega_\mathcal{B},
	$$
where $\pi_{\MM_H}$ denotes the quotient map of $\MM_H$ and $\iota$ is the inclusion into $X$ (see \cite{MW}, p. 123). In fact, $\MM_H/H$ is even a Kähler manifold. The key step is to verify that the almost complex structure naturally arising from the construction of the symplectic reduction of $X$ is integrable (see \cite{GS}, p. 522). This follows as a particular case of a more general result in \cite{HS}, which states that for any closed Lie subgroup $H$ of a Kähler manifold's holomorphic isometry group acting in a Hamiltonian fashion, the symplectic reduction $\MM_H/H$ is a Kähler space. We are interested in whether the symplectic reduction of $X$ is a Stein manifold and address this question for two particular classes of subgroups in \textbf{Section \ref{section3}}.\vspace{6pt}
	
First, we consider a homogeneous Siegel domain $$T_\Omega:=V+i\Omega\subset V^\CC:=V+iV$$ of the first kind, where $\Omega \subset V$ is a homogeneous proper open convex cone in a finite-dimensional Euclidean vector space $V$. Further, let $H$ be a connected subgroup of the translation group $V$, i.e., a vector subspace of $V$ and $$Z:=H^\CC+T_\Omega=V+i(H+\Omega).$$ We observe that $Z$ is an open convex subset of $V^\CC$ and $Z/H^\CC$ is open in the quotient vector space $V^\CC/H^\CC$. Let $\mu_H:T_\Omega\longrightarrow\mathfrak{h}^*$ be the momentum map (see Sect. \ref{twofour} and \ref{Example1}). The inclusion $\iota$ of $\MM_H$ into $Z$ induces a map $\overline\iota:\MM_{H}/H\longrightarrow Z/H^\CC$ and we have the following commutative diagram:\vspace{6pt}
\begin{center}
	\begin{tikzcd}
		\MM_H \arrow[r, hook, "\iota"] \arrow[d]
		& Z \arrow[d] \\
		\MM_H/H \arrow[r,  "\overline\iota"]
		& Z/H^\CC
	\end{tikzcd}
\end{center} As a first result we prove the following theorem.
	\begin{THEO}\label{THEO:theoone}(Thm. \ref{Theo:inclusiontrans})
		Let $T_\Omega$ be a homogeneous Siegel domain of the first kind and $H<V$ a connected subgroup of the translation group. Set $Z:=H^\CC+T_\Omega$. Then the map $\overline\iota:\MM_H/H\longrightarrow Z/H^\CC$ induced by the inclusion $\iota:\MM_{H}\longrightarrow Z$ is biholomorphic. In particular, the symplectic reduction $\MM_H/H$ is a Stein manifold biholomorphic to a Siegel domain of the first kind.
	\end{THEO}

Next, we consider a general homogeneous Siegel domain $X$ in $\CC^d$. Let $B$ be a maximal simply connected split solvable subgroup of $\textup{Aut}_{\mathcal{O}}(X)$ and let $H$ be a closed subgroup of $B$ with Lie algebra $\mathfrak{h}$. For all $x\in \MM_H$, we set $W_{x}:=T_{x}\MM_H\cap JT_{x}\MM_H$, where $J$ denotes the complex structure of $X$. Further, for $x\in X$ let 
	$$
	\mathfrak{h}\cdot x:=\{\xi_X(x)\ |\ \xi\in\mathfrak{h}\}=T_x(H\cdot x)\subset T_xX
	$$ 
be the tangent space of the $H$-orbit at $x$. For $x\in \MM_H$, we have  
	$$
	\textup{Ker}(d\mu_H(x))=T_{x}\MM_H=(\mathfrak{h}\cdot x)^{\perp_{\omega_\mathcal{B}}}=\mathfrak{h}\cdot x\oplus W_{x},
	$$
where $\perp_{\omega_\mathcal{B}}$ denotes perpendicular with respect to $\omega_\mathcal{B}$. 
	\begin{Deff}
		We say that the group $H$ satisfies the Lie condition if
		\begin{enumerate}[i)]
			\item $\mathcal{M}_H$ is connected,
			\item there exists $x_0$ in $\MM_H$ and a connected subgroup $S=\exp(\mathfrak{s})$ of $B$ such that $S\cdot x_0\subset\MM_H$ and $\mathfrak{s}\cdot x_0=W_{x_0}$.
		\end{enumerate}
	\end{Deff}
As $B$ is split solvable, $S$ is necessarily closed (see \cite{VinLie}, p. 52). In fact, we prove that the $S$-orbit through $x_0$ is a closed complex submanifold of $X$. Our main result is the following.
	\begin{THEO}\label{THEO:theotwo}(Thm. \ref{Theo:Quotient})
		Let $X\subset\CC^d$ be a homogeneous Siegel domain and $H$ a closed subgroup of $B$ which satisfies the Lie condition. Then
		\begin{enumerate}[i)]
			\item $\MM_H=HS\cdot x_0$,
			\item the restriction of the quotient map $\pi:\MM_{H}\longrightarrow\MM_H/H$ to $S\cdot x_0$ is a biholomorphic Kähler isometry.
		\end{enumerate}
		In particular, the symplectic reduction $\MM_H/H$ is Stein.
	\end{THEO}

\subsection*{Outline of the paper}
In \textbf{Section \ref{section2}} we recall some facts about tube domains and their relation with homogeneous bounded domains in $\CC^d$. Further, we introduce the momentum map for tube domains and compute a relevant example. In \textbf{Section \ref{threeone}} we state and prove Theorem \ref{THEO:theoone} and in \textbf{Section \ref{threethree}} we prove Theorem \ref{THEO:theotwo}.

\section{Hamiltonian geometry of tube domains}\label{section2}
In this section we analyse some properties of open convex cones and discuss the Hamiltonian geometry of tube domains. Then we recall those results on homogeneous bounded domains which are needed in the sequel. An explicit computation is pointed out in the case of the Lorentz tube. Details can be found in \cite{P}, \cite{Kan} and \cite{F}.
\subsection{Open convex cones}\label{twoone}
Let $V$ be a finite-dimensional Euclidean vector space. A subset $ C\subset V$ is called a cone if $x\in C$ and $\lambda\in\RR^{>0}$ imply $\lambda x\in C$. It is called convex if for $x,y\in C$ and $0<\lambda<1$ we have $\lambda x+(1-\lambda)y\in C$. It is easy to show that a subset $C$ is a convex cone if and only if $x,y\in C$ and $\lambda,\mu\in\RR^{>0}$ imply $\lambda x+\mu y\in C$. The definition of a cone implies that the interior of a cone $C$ is not empty if and only if $C$ contains a basis of $V$. In this case the interior is itself an open cone. For any cone $C$ the topological closure $\overline C$ is a cone.
	\begin{Def}
		Let $\Omega\subset V$ be an open cone.
		\begin{enumerate}[i)]
			\item We call $\Omega$ proper if $\overline\Omega\cap-\overline\Omega=\{0\}$.
			\item For an inner product $g:V\times V\rightarrow\RR$ the open dual cone of $\Omega$ is defined by
			\begin{align*}
				\Omega^*:=\{v\in V\ |\ g(v,\omega)>0\ \text{for all }\omega\in\overline\Omega\setminus\{0\}\}.
			\end{align*} 
			\item We call $\Omega$ self-dual if there exists an inner product $g$ such that $\Omega=\Omega^*$. 
			\item The linear automorphism group of $\Omega$ is defined by
			\begin{align*}
				G_\Omega:=\{g\in\textup{GL}(V)\ |\ g\cdot\Omega=\Omega\}.
			\end{align*}
			We call $\Omega$ homogeneous if $G_\Omega$ acts transitively on it.
			\item We call $\Omega$ symmetric if it is self-dual and homogeneous.
		\end{enumerate}
	\end{Def}
	
	\begin{Rk}\textit{(Properties of the dual cone.)}\mbox{}
		\begin{enumerate}[i)]
			\item The dual cone varies depending on the choice of the inner product.
			\item The dual cone of $\overline\Omega$ is by definition the cone
			\begin{align*}
				\overline\Omega^*=\{v\in V\ |\ g(v,\omega)\geq0\ \text{for all }\omega\in\overline\Omega\}.
			\end{align*} 
			\item For an open cone $\Omega$ the dual cone $\Omega^*$ is open. It is also the interior of the dual cone of $\overline\Omega$. Moreover, by the bilinearity of $g$, the dual cone is convex.
			\item A self-dual open cone is always convex and proper (see \cite{F}, p. 4).
		\end{enumerate}
	\end{Rk}
	\begin{Ex}\label{Ex:cones} On $\RR^{d+1}$ an inner product is given by the Euclidean scalar product
			\begin{align*}
				\left\langle v,w\right\rangle:=v_0w_0+v_1w_1+\dots+v_dw_d,\ v,w\in\RR^{d+1}.
			\end{align*}
			We introduce the Lorentz product
			\begin{align*}	
				\left\langle v,w\right\rangle_{1,d}:=v_0w_0-v_1w_1-\dots-v_dw_d,\ v,w\in\RR^{d+1}.
			\end{align*}
			The domain
			\begin{align*}
				\Omega=\{\omega\in\RR^{d+1}\ |\ \omega_0>0,\ \left\langle \omega,\omega\right\rangle_{1,d}>0\}
			\end{align*}
			is called the Lorentz cone and is a symmetric cone in $\RR^{d+1}$ (see \cite{F}, pp.~7-8).
	\end{Ex}
The orthogonal complement of a vector subspace $H$ in $V$ is defined by
	$$
	H^{\perp_g}:=\left\{v\in V\ |\ g(v,h)=0\textup{ for all }h\in H\right\}.
	$$ 
	
	\begin{Theo}\label{Theo:compact}
		Let $\Omega\subset V$ be a proper open convex cone and let $H\subset V$ be a vector subspace such that $H^{\perp_g}\cap\Omega^*\neq\emptyset$. If $K\subset\Omega$ is compact, the set 
		$$
		(H+K)\cap\overline\Omega=\left\{h+k\in V\ |\ h\in H,\ k\in K\right\}\cap\overline\Omega
		$$ 
		is compact.
	\end{Theo}
In particular, for $K=\{\omega\}$ the intersection $(H+\omega)\cap\overline\Omega$ is compact. For the proof we need the following two lemmata (for the proof of the first lemma see \cite{F}, p. 3).
	\begin{Lem}\label{Lem:constant}
		Let $U\subset\Omega^*$ be a compact subset. Then there exists a constant $p>0$ such that for all $\omega\in\overline\Omega$ and for all $y\in U$, we have $p||\omega||\leq g(\omega,y)$, where $||\omega||:=\sqrt{g(\omega,\omega)}$ is the induced norm on V.
	\end{Lem}
	\begin{Lem}\label{Lem:Coneperp}
		Let $\Omega\subset V$ be a proper open convex cone and $H\subset V$ a vector subspace. We have 
		$
		H^{\perp_g}\cap\Omega^*\neq\emptyset \Leftrightarrow H\cap\overline\Omega=\{0\}.
		$
	\end{Lem}
\Proof 
i) Assume that $H^{\perp_g}\cap\Omega^*\neq\emptyset$. Let $\omega\in H^{\perp_g}\cap\Omega^*$ and assume that there exists $h\in H\cap\overline\Omega$ with $h\neq 0$. Then $g(\omega,h)>0$, giving a contradiction.\vspace{6pt}
	
\noindent ii) Assume that $H\cap\overline\Omega=\{0\}$. This implies that $H\cap\left(\overline{\Omega}\setminus\{0\}\right)=\emptyset$. Since $\overline{\Omega}\setminus\{0\}$ is convex, there exists a $\lambda\in V^*$ such that $H\subset\textup{Ker}(\lambda)$ and $\lambda(\omega)>0$ for all $\omega\in\overline\Omega\setminus\{0\}$ (see \cite{R}, §11). In particular, for the gradient $x:=\nabla\lambda$ which is defined with respect to the inner product $g$, we have $x\in\Omega^*$, hence $H^{\perp_g}\cap\Omega^*\neq\emptyset$.
\qed\vspace{6pt}
	
\textit{Proof of Theorem \ref{Theo:compact}.} 
Let $\Omega\subset V$ be a proper open convex cone, $H\subset V$ a vector subspace such that $H^{\perp_g}\cap\Omega^*\neq\emptyset$ and $K\subset\Omega$ a compact subset. Let $\omega\in K$ and $y\in H^{\perp_g}\cap\Omega^*$. Lemma \ref{Lem:Coneperp} implies that for all $h\in H$ with $h+\omega\in\overline\Omega$, we have $h+\omega\neq 0$. For the compact set $U=\left\{y\right\}$, Lemma \ref{Lem:constant} implies that there exists a constant $p>0$ such that $p||h+\omega||\leq g(h+\omega,y)$, where $p$ does not depend on $h+\omega$. As $y\in H^{\perp_g}$, one has $||h+\omega||\leq\frac{1}{p}g(\omega,y)$. The compact set $K$ is bounded by a constant $r>0$, so the Cauchy-Schwarz inequality implies
	\begin{align*}
		||h+\omega||\leq\frac{1}{p}g(\omega,y)\leq\frac{1}{p}||\omega||\ ||y||\leq\frac{r}{p}||y||.
	\end{align*}
This shows that the set $(H+K)\cap\overline\Omega$ is bounded. Moreover, $H+K$ is a sum of a closed and a compact set. Therefore $(H+K)\cap\overline\Omega$ is closed, thus compact.
\qed
	
\subsection{The characteristic function of a proper open convex cone}\label{twotwo}
Let $V$ be a finite-dimensional vector space endowed with an inner product $g$, i.e., a Euclidean vector space and let $\Omega\subset V$ be a proper open convex cone. The characteristic function of $\Omega$ is the map
	\begin{align*}
		\varphi:\Omega\longrightarrow\RR,\ \varphi(\omega)=\int\limits_{\Omega^*}e^{-g(\omega,y)}dy,
	\end{align*}
where $dy$ denotes the Euclidean measure on $V$. The following proposition summarizes the main properties of $\varphi$ (see \cite{F}, pp. 10-11).
	\begin{Prop}\label{Prop:propertiesphi}\mbox{}
		\begin{enumerate}[i)]
			\item The characteristic function $\varphi$ is real analytic. 
			\item For every sequence $(\omega_n)_{n\in\NN}\subset\Omega$ converging to a boundary point $\omega_0\in\partial\Omega$, we have $\lim\limits_{n\rightarrow\infty}\varphi(\omega_n)=\infty.$
			\item For all $g\in G_\Omega$, we have $\varphi(g\cdot \omega)=|\det g|^{-1}\varphi(\omega)$.
			\item The function $\log(\varphi)$ is strictly convex.
		\end{enumerate}
	\end{Prop}
In particular, $\varphi$ is invariant under any subgroup of $G_\Omega\cap\textup{SL}(V)$ and homogeneous of degree $-\dim_{\RR} V$. The gradient of $\log(\varphi)$ with respect to the inner product $g$ defines the dual map
	\begin{align*}
		\psi:\Omega\longrightarrow V,\ \psi(\omega)=-\nabla\log(\varphi(\omega))=:\omega^*.
	\end{align*}
One has the following proposition (see \cite{F}, pp. 13-14).
	\begin{Prop}\label{Prop:properties}\mbox{}
		\begin{enumerate}[i)]
			\item For all $\omega\in\Omega$, we have $g(\omega,\omega^*)=\dim_{\RR} V$.
			\item The dual map $\psi$ is a bijection of $\Omega$ onto $\Omega^*$.
		\end{enumerate}
	\end{Prop}
	\begin{Ex}\label{Ex:conedual}
		Let $\Omega\subset{\RR^{d+1}}$ be the Lorentz cone (see Ex. \ref{Ex:cones}). The dual mapping is given by
		\begin{align*}
			(\omega_0,\dots,\omega_d)^{*}=\frac{d+1}{\left\langle \omega,\omega\right\rangle_{1,d}}(\omega_0,-\omega_1,\dots,-\omega_d)\ \text{(see \cite{F}, p. 15)}.
		\end{align*}
	\end{Ex}
\subsection{Homogeneous bounded domains and Siegel domains}\label{twothree}
Let $X\subset\CC^d$ be a homogeneous bounded domain. Our main interest is in the Hamiltonian geometry of such domains. For this we rely on the fact that they are biholomorphic to Siegel domains of the first or second kind (see \cite{Kan} and \cite{P}). We give a precise description of Siegel domains of the first kind. For more details on Siegel domains of the second kind see \cite{P}, Ch. 1, Sect. 2.
	\begin{Def}
		Let $V$ be a finite-dimensional Euclidean vector space and let $\Omega\subset V$ be a proper open convex cone. Further, let $V^\CC:=V\oplus iV$. The tube domain 
		$$
		T_\Omega:=V+i\Omega\subset V^\CC
		$$ is called Siegel domain of the first kind.
	\end{Def}
Every Siegel domain of the first kind is biholomorphic to a bounded domain (see \cite{P}, p. 16). Vinberg, Gindikin and Pyatetskii-Shapiro have proved the following (see \cite{P}, Appendix).
	\begin{Theo}
		Every homogeneous bounded domain $X\subset\CC^d$ is biholomorphic to a Siegel domain of the first or second kind.
	\end{Theo}
Let $\textup{Aut}_{\mathcal{O}}(X)$ be the group of holomorphic automorphisms of $X$ and $G$ its connected component containing the identity. For a fixed base point $x_0\in X$ we have a decomposition $G=BG_{x_0}$ (see \cite{Kan} and \cite{Vin}). The split solvable subgroup $B$ is isomorphic to the semi-direct product $A\ltimes N$, where $N$ denotes the nilradical of $B$ and $A\cong (R^{>0})^d$. The group $B$ acts on the Siegel domain realization of $X$ by affine transformations (see \cite{Kan}, Ch.~1) and this action extends to all of $\CC^d$. Moreover, the $N$-action extends to an algebraic action of the universal complexification $N^\CC$ on $\CC^d$ (see \cite{Kan}, Ch. 2 and 3).
	\begin{Rk}\label{Rk:connected}
		Any connected Lie subgroup of $B$ is a closed subgroup (see \cite{VinLie}, p. 52).
	\end{Rk}
Let $T_\Omega=V+i\Omega$ be a homogeneous Siegel domain of the first kind. The subgroup $\textup{Aff}(T_\Omega)$ of $G$ consisting of affine transformations of $T_\Omega$ is given by the semi-direct product of the linear automorphism group of the cone $G_\Omega$ and the group of real translations $V$. For a fixed base point $x_0=0+i\omega_0\in T_\Omega$ we have a decomposition $G=B(G_\Omega)_{\omega_0}$ (see \cite{Kan}, Ch. 1). 
	\begin{Ex}\label{Ex:lorentztube}
		Consider the Lorentz cone $\Omega\subset\RR^{d+1}$ (see Ex. \ref{Ex:cones}). The Lorentz tube is the tube domain $T_\Omega=\RR^{d+1}+i\Omega$. The Lie ball 
		$$
		X=\left\{z\in\CC^{d+1}\ |\ |z|^2+\sqrt{|z|^4-|z_0^2+\dots+z_d^2|^2}<1\right\}\subset\CC^{d+1}
		$$
		is biholomorphic to $T_\Omega$. The connected component of the identity in the automorphism group is isomorphic to $G=\textup{SO}_0(2,{d+1})$ and the subgroup of affine transformations is $G_\Omega\ltimes\RR^{d+1}$ with $G_\Omega=\RR^{>0}\times\textup{SO}_0(1,d)$ (see \cite{F}, Ch. X.5). The group $\textup{SO}_0(1,d)$ is semi-simple. We consider an Iwasawa decomposition $\textup{SO}_0(1,d)=K_\Omega A_\Omega N_\Omega$, where $A_\Omega\cong\RR^{>0}$ is abelian and $N_\Omega$ is unipotent (see \cite{Kn}, Ch. VI.4). In particular, $B=\left(\RR^{>0}\times\RR^{>0}\right)\ltimes N$ with $N=N_\Omega\ltimes\RR^{d+1}$. 
	\end{Ex}

\subsection{Hamiltonian geometry of tube domains}\label{twofour}\emergencystretch=1em
The Bergman metric on $T_\Omega$ defines a Kähler form $\omega_\mathcal{B}$ on $T_\Omega$. The Bergman kernel of $T_\Omega$ (see \cite{F}, p. 177) is up to a constant given by 
	\begin{align*}
		K(x,x)=c\varphi(2\im x)^2,\ c>0,\ x\in T_\Omega.
	\end{align*}
Set $\textup{SL}_\Omega:=G_\Omega\cap\textup{SL}(V)$. We observe that by Proposition \ref{Prop:propertiesphi} the kernel is $\left(\textup{SL}_\Omega\ltimes V\right)$-invariant and defines an $\left(\textup{SL}_\Omega\ltimes V\right)$-invariant Kähler potential
	\begin{align*}
		\rho:T_\Omega\longrightarrow\RR,\ \rho(x)=c_B\log(\varphi(\im x))+\widetilde c_B
	\end{align*}
with respect to the Kähler form $\omega_\mathcal{B}$. That is, $\omega_\mathcal{B}=-dd^c\rho$, where $d^c\rho:=d\rho\circ J$ and $J$ is the complex structure on $T_\Omega$. Let $\mathfrak{sl}_\Omega\oplus V$ be the Lie algebra of $\textup{SL}_\Omega\ltimes V$. The potential $\rho$ induces a momentum map with respect to the $\left(\textup{SL}_\Omega\ltimes V\right)$-action on $T_\Omega$, given by
	\begin{align*}
		\mu^\xi:T_\Omega\longrightarrow\RR,\  \mu^\xi(x)=d^c\rho(\xi_{T_{\Omega}}(x)),\  \xi\in\mathfrak{sl}_\Omega\oplus V.
	\end{align*}
We want to describe the momentum map more precisely in the case where $c_B=1$ and $\widetilde c_B=0$, i.e., $\rho(x)=\log(\varphi(\im x))$. For $x=v+i\omega\in T_\Omega$ and $\xi=\xi_\Omega+\xi_V\in\mathfrak{sl}_\Omega\oplus V$, we compute
	\begin{align*}
		\mu^\xi(x)&=d^c\rho(\xi_{T_{\Omega}}(x))=\frac{d}{dt}\Big|_{t=0}\rho(\exp(it\xi)\cdot x)=\frac{d}{dt}\Big|_{t=0}\rho(it\xi_V+\exp(it\xi_\Omega)\cdot x)\\
		&=\frac{d}{dt}\Big|_{t=0}\log\left(\varphi\left(t\xi_V+\im\left(\exp(it\xi_\Omega)\cdot x\right)\right)\right)\\
		&=g(\nabla\log(\varphi(\omega)),\xi_V)+g(\nabla\log(\varphi(\omega)),\im\left(i\xi_\Omega\cdot x\right))\\
		&=g(\nabla\log(\varphi(\omega)),\xi_V)+g(\nabla\log(\varphi(\omega)),\xi_\Omega\cdot v)\\
		&=-g(\omega^*,\xi_V+\xi_\Omega\cdot v).
	\end{align*}
This computation and the fact that $\textup{d}\mu^\xi(x)=\omega_\mathcal{B}(\xi_{T_\Omega}(x),\cdot)$ imply the following.
	\begin{Prop}
		Let $\mathfrak{g}_\Omega\oplus V$ be the Lie algebra of $\textup{Aff}(T_\Omega)=G_\Omega\ltimes V$ and $\xi=\xi_\Omega+\xi_V\in\mathfrak{g}_\Omega\oplus V$. The map 
		$$
		\mu^\xi:T_\Omega\longrightarrow\RR,\ \mu^\xi(v+i\omega)=-g(\omega^*,\xi_V+\xi_\Omega\cdot v)
		$$
		defines a momentum map for the ${\textup{Aff}(T_\Omega)}$-action on $T_\Omega$ with respect to the Kähler form $\omega_\mathcal{B}$.
	\end{Prop}
For any Lie subgroup $H$ of $\textup{Aff}(T_\Omega)$, the inclusion of the Lie algebra of $H$ into $\mathfrak{g}_\Omega\oplus V$ induces a momentum map $\mu_H$ by restriction. In the next section we analyse the zero level set $\MM_H:=\mu^{-1}_{H}(0)$ of $\mu_H$ for an important class of subgroups.
	
\subsection{Example: Momentum map for subgroups of the translation group}\label{Example1}
Let $H<V$ be a connected subgroup, i.e., a vector subspace of $V$. By the previous computations the momentum map for the $H$-action is given by
	$$
	\mu^\xi_{H}(v+i\omega)=-g(\omega^*,\xi).
	$$
The zero level set of $\mu_H$ is given by $\MM_{H}=V+iC_{H}$, where $C_{H}:=\psi^{-1}(H^{\perp_g}\cap\Omega^*)$. Here $\psi:\Omega\longrightarrow\Omega^*$ denotes the dual mapping that we introduced in the Section \ref{twotwo}. We observe that 
	$$
	\MM_{H}\neq\emptyset\Leftrightarrow C_{H}\neq\emptyset\Leftrightarrow H^{\perp_g}\cap\Omega^*\neq\emptyset.
	$$ 
\newpage
	\begin{Prop}\label{Prop:empty}\mbox{}
		\begin{enumerate}[i)]
			\item The set $C_H$ is a connected cone. In particular, $\MM_H$ is a connected set.
			\item $C_{H}\neq\emptyset\Leftrightarrow H\cap\overline\Omega=\{0\}$. 
		\end{enumerate}
	\end{Prop}
\Proof 
 It follows from Proposition \ref{Prop:properties} (i) that $(\lambda\omega)^*=\frac{1}{\lambda}\omega^*$ for every real positive $\lambda$. This implies that $C_H$ is a cone. The inverse $\psi^{-1}$ of the dual map is continuous and the set $H^{\perp_g}\cap\Omega^*$ is connected since $\Omega^*$ is convex. Thus, $C_H$ is connected and so is $\MM_H$. The second statement follows from Lemma~\ref{Lem:Coneperp}.
\qed\vspace{6pt}
	\begin{Ex}
		Consider the Lorentz tube $T_\Omega=\RR^{d+1}+i\Omega$ (see Ex. \ref{Ex:cones} and \ref{Ex:conedual}) and a vector subspace $H\subset\RR^{d+1}$ with $H^{\perp_{\langle\ \rangle_{1,d}}}\cap\Omega\neq\emptyset$. For $\omega\in\Omega$ and $\xi\in H$, we have
		\begin{align*}
			\left\langle\omega^*,\xi\right\rangle=\frac{d+1}{\left\langle \omega,\omega\right\rangle_{1,d}}\left\langle \omega,\xi\right\rangle_{1,d}
		\end{align*}
		so the momentum map for the $H$-action on $T_\Omega$ is given by
		\begin{align*}
			\mu_{H}^\xi(v+i\omega)=-\frac{d+1}{\left\langle \omega,\omega\right\rangle_{1,d}}\left\langle\omega,\xi\right\rangle_{1,d}
		\end{align*}
		and therefore
		\begin{align*}
			\MM_{H}=\RR^{d+1}+i(H^{\perp_{\left\langle\ \right\rangle_{1,d}}}\cap\Omega).
		\end{align*}
	\end{Ex}
In the sequel we assume that $C_{H}\neq\emptyset$. We have already seen in Theorem~\ref{Theo:compact} that for all $\omega\in\Omega$, the intersection of $\overline\Omega$ with $H+\omega$ is compact. For the intersection of $C_{H}$ with such an affine space we get a much stronger result.
	\begin{Theo}\label{Theo:section}
		For all $\omega\in C_{H}$, we have $(H+\omega)\cap C_{H}=\{\omega\}$.
	\end{Theo}	
\Proof
Let $\omega\in C_{H}$, i.e., $\omega^*\in H^{\perp_g}\cap\Omega^*$. Then $g(\omega,\omega^*)=d$, where $d=\dim_\RR V$ (see Prop. \ref{Prop:properties} (i)). We define the sets
	\begin{align*}
		H_d(\omega^*):=\{v\in V\ |\ g(v,\omega^*)=d\}\text{ and }Q(\omega^*):=H_d(\omega^*)\cap\Omega.
	\end{align*}
Then $\omega\in Q(\omega^*)$ is the unique minimum of $\varphi|_{Q(\omega^*)}$ (see \cite{F}, pp. 13-14). Since $g(\omega^*,h)=0$ for all $h\in H$, it follows that $(H+\omega)\subset H_d(\omega^*)$ and that $\omega$ is the unique minimum of $\varphi|_{(H+\omega)\cap\Omega}$. If $z=h+\omega\in (H+\omega)\cap C_{H}$, then $z^*\in H^{\perp_g}\cap\Omega^*$. Further, $z$ is the unique minimum of $\varphi|_{Q(z^*)}$. But $z^*\in H^{\perp_g}$ implies $d=g(\widetilde h+\omega,z^*)=g(\omega,z^*)$ for all $\widetilde h\in H$. This shows that $H+\omega\subset H_d(z^*)$. In conclusion $z$ has to be the unique minimum of $\varphi$ on $(H+\omega)\cap\Omega\subset Q(z^*)$ which shows that $z=\omega$ and therefore $h=0$.
\qed\vspace{6pt}
	
Let $H^\CC<V^\CC$ be the complexification of $H$. The group $H^\CC$ acts on $V^\CC$ by translations. 
	\begin{Coro}\label{Coro:inject}\mbox{}
		\begin{enumerate}[i)]
			\item For all $x\in \MM_{H}$, we have $H^\CC\cdot x\cap\MM_{H}=H\cdot x$.
			\item The map $\overline\iota:\MM_{H}/H\longrightarrow V^\CC/H^\CC$ induced by the inclusion $\iota:\MM_H\longrightarrow V^\CC$ is a biholomorphism onto its open image.
		\end{enumerate}
	\end{Coro} 
\Proof 
i) Let $x=v+i\omega\in\MM_{H}$ and $h_1+ih_2\in H^\CC$. Then 
$$
(h_1+ih_2)\cdot x=(h_1+v)+i(h_2+\omega)
$$
belongs to $\MM_{H}$ if and only if $h_2+\omega\in C_{H}$ and Theorem \ref{Theo:section} implies $h_2=0$, i.e., $h_1+ih_2=h_1\in H$.\vspace{6pt}

ii) Since $H$ acts freely on $T_\Omega$, the momentum map $\mu_H$ is a submersion and we have
    $$
	\dim_\RR V^\CC=\dim_\RR T_\Omega=\dim_\RR\MM_{H}+\dim_\RR H. 
	$$
This and (i) imply that the projection of $\MM_H$ to $V^\CC/H^\CC$ factorizes to an injective local diffeomorphism which induces the complex structure of the symplectic reduction $\MM_H/H$ from that of $V^\CC/H^\CC$ (cf. \cite{Kur}, Sect. 2).
\qed
	
\section{Symplectic reduction for two particular classes of subgroups}\label{section3}
Let $X\subset\CC^d$ be a homogeneous Siegel domain endowed with the Kähler form $\omega_\mathcal{B}$ induced by the Bergman metric. The group of holomorphic automorphisms $\textup{Aut}_{\mathcal{O}}(X)$ acts in a Hamiltonian fashion on $X$ with momentum map $\mu$. Given a closed subgroup $H$ of $\textup{Aut}_{\mathcal{O}}(X)$, consider the restricted momentum map $\mu_H$ and the associated symplectic reduction $\MM_H/H$. The aim of this section is to show that the symplectic reduction of $X$ is a Stein manifold when $H$ belongs to two different classes of closed subgroups of $\textup{Aut}_{\mathcal{O}}(X)$. 
	
\subsection{Symplectic reduction for subgroups of the translation group}\label{threeone}
Let $T_\Omega:=V+i\Omega$ be a homogeneous Siegel domain of the first kind. Let $H<V$ be a connected subgroup of the translation group, i.e., a vector subspace of $V$ and let $H^\CC<V^\CC$ be its complexification. Set $Z:=H^\CC+T_\Omega$. In this section we prove that the symplectic reduction $\MM_H/H$ is biholomorphic to $Z/H^\CC$. Its Steiness will follow by showing that $Z/H^\CC$ is biholomorphic to a Siegel domain of the first kind.\\
In the sequel we assume that $\MM_{H}\neq\emptyset$, i.e., $H\cap\overline\Omega=\{0\}$ (see Prop. \ref{Prop:empty}) and extend $g$ to a Hermitian inner product $h$ on $V^\CC$.
	\begin{Prop}\label{Prop:Siegel}
		Let $T_\Omega$ be a homogeneous Siegel domain of the first kind and let $H<V$ be a connected subgroup.
		\begin{enumerate}[i)]
			\item The quotient $Z/H^\CC$ is Stein and $Z=H^\CC\oplus S$, where $S$ is the orthogonal projection of $Z$ to $H^{\CC^{\perp_h}}$. 
			\item $S=H^{\perp_g}+i\left(H^{\perp_g}\cap(H+\Omega)\right)$.
			\item $S$ is a Siegel domain of the first kind.
		\end{enumerate}
	\end{Prop}
\Proof
i) Let $P:V^\CC\longrightarrow H^{\CC^{\perp_h}}$ be the orthogonal projection. Note that $Z=H^\CC+T_\Omega=V+i(H+\Omega)$ is open and convex in $V^\CC$. Therefore, $S:=P(Z)$ is open and convex in $H^{\CC^{\perp_h}}$. It follows that $Z=H^\CC\oplus S$ is also convex, implying that $Z$ is Stein. \vspace{6pt}
	
ii) Let $v+i\omega\in T_\Omega$ and $z=h_1+v+i(h_2+\omega)\in Z$, where $h_1,h_2\in H$. We have $V=H\oplus H^{\perp_g}$. Accordingly, $\omega$ admits a decomposition $\omega=\omega_H+\omega_{H^\perp}$ with $\omega_H\in H$ and $\omega_{H^\perp}\in H^{\perp_g}\cap(H+\Omega)$. Consequently, we have
	$$
	h_2+\omega=(h_2+\omega_H)+\omega_{H^\perp}.
	$$
Therefore, we get
	\begin{align*}
		P(z)&=h^\perp_{h_1+v}+i\omega_{H^\perp}\in H^{\perp_g}+i\left(H^{\perp_g}\cap(H+\Omega)\right).
	\end{align*}
	
iii) We have to show that $S$ is a tube domain over a proper open convex cone. Since $H+\Omega$ is a convex cone, $H^{\perp_g}\cap(H+\Omega)$ is an open convex cone in $H^{\perp_g}$. It remains to show that $H^{\perp_g}\cap(H+\Omega)$ is proper, that is, 
	\begin{align*}
		\left(\overline{H^{\perp_g}\cap(H+\Omega)}\right)\cap-\left(\overline{H^{\perp_g}\cap(H+\Omega)}\right)=\{0\}.
	\end{align*}
First, we show that 	
	\begin{equation}\tag{1}\label{eq:one}
		\overline{H^{\perp_g}\cap(H+\Omega)}=H^{\perp_g}\cap\overline{(H+\Omega)}.
	\end{equation}
Let $(x_n)_{n\in\NN}\subset H^{\perp_g}\cap(H+\Omega)$ be a convergent sequence. The inner product $g$ is continuous, hence 
	$$
	g(\lim\limits_{n\to\infty}x_n,\widetilde h)=\lim\limits_{n\to\infty}g(x_n,\widetilde h)=0
	$$
for all $\widetilde h\in H$ which proves that $x=\lim\limits_{n\to\infty}x_n\in H^{\perp_g}$. Hence, the desired equality follows.\\
Next, we show that 
	\begin{equation}\tag{2}\label{eq:two}
		\overline{H+\Omega}=H+\overline\Omega.
    \end{equation}
This will imply that
	\begin{equation}\tag{3}\label{eq:three}
		\overline{H^{\perp_g}\cap(H+\Omega)}\overset{(\ref{eq:one})}{=}H^{\perp_g}\cap\overline{(H+\Omega)}\overset{(\ref{eq:two})}{=}H^{\perp_g}\cap(H+\overline\Omega).
    \end{equation}
Let 
	$$
	(x_n)_{n\in\NN}=(h_n+\omega_n)_{n\in\NN}\subset H+\overline\Omega\setminus\{0\}
	$$ 
be a convergent sequence with $\lim\limits_{n\to\infty}x_n=x\neq 0$. Then there exist sequences $(\lambda_n)_{n\in\NN}\subset\RR^{>0}$ and $(k_n)_{n\in\NN}\subset S(V)\cap(\overline\Omega\setminus\{0\})$ such that $\omega_n=\lambda_n\cdot k_n$ for all $n\in\NN$, where $S(V)$ denotes the unit sphere in $V$. Since $S(V)\cap(\overline\Omega\setminus\{0\})$ is compact, we can assume that $(k_n)_{n\in\NN}$ is convergent. In particular, $$
(x_n,k_n)_{n\in\NN}=\left((\lambda_n,h_n)\cdot k_n,k_n\right)_{n\in\NN}
	$$ 
is a convergent sequence in $V\setminus\{0\}\times V\setminus\{0\}$. Since $\RR^{>0}$ acts properly on $V\setminus\{0\}$ and $H$ acts properly on $V$, $(\lambda_n,h_n)_{n\in\NN}$ has a convergent subsequence. In particular, $\omega_n=\lambda_n\cdot k_n$ has a convergent subsequence. This implies that $x\in(H+\overline\Omega)\setminus\{0\}$. Summarized, we have 
	$$
	\overline{H+\overline{\Omega}\setminus\{0\}}=H+\overline\Omega
	$$
and therefore $\overline{H+\Omega}\subset H+\overline\Omega$. Furthermore, we have $H+\overline\Omega\subset\overline{H+\Omega}$ which proves the equality.\\ 
Finally, we show that 
	\begin{equation}\tag{4}\label{eq:four}
		(H+\overline\Omega)\cap-(H+\overline\Omega)=H.
    \end{equation}
Let $h+\omega\in (H+\overline\Omega)\cap-(H+\overline\Omega)$. Then there exists an $\widetilde h+\widetilde\omega$ such that $h+\omega=-(\widetilde h+\widetilde\omega)$. This implies that $\omega+\widetilde\omega\in H\cap\overline\Omega$. By Proposition \ref{Prop:empty}, $H\cap\overline\Omega=\{0\}$, hence $\omega=-\widetilde\omega$. Since $\Omega$ is proper, this implies $\omega=0$.\\ 
Putting everything together, we obtain
	\begin{align*}
		\left(\overline{H^{\perp_g}\cap(H+\Omega)}\right)\cap-\left(\overline{H^{\perp_g}\cap(H+\Omega)}\right)&\overset{(\ref{eq:three})}{=}H^{\perp_g}\cap\left((H+\overline\Omega)\cap-(H+\overline\Omega)\right)\\
		&\overset{(\ref{eq:four})}{=}H^{\perp_g}\cap H=\{0\}.
	\end{align*}
\qed\vspace{6pt}	
	
Now it remains to show that $\MM_H/H$ is biholomorphic to $Z/H^\CC$. First, we prove the following.
	\begin{Prop}
		The set $\left\{h\in H^\CC\ |\ h\cdot x\in T_\Omega\right\}$ is connected for every $x\in T_\Omega$.
	\end{Prop}
\Proof 
Let $x=v+i\omega\in T_\Omega$ and $h=h_1+ih_2\in H^\CC$. We observe that 
	$$
	h\cdot x=h_1+v+i(h_2+\omega)\in T_\Omega\Leftrightarrow h_2+\omega\in\Omega.
	$$
Therefore,
	$$
	\left\{h\in H^\CC\ |\ h\cdot x\in T_\Omega\right\}=\left\{h_1+ih_2\in H^\CC\ |\ h_2+\omega\in\Omega\right\}
	$$
and the latter set is connected since $(H+\omega)\cap\Omega$ is connected.
\qed\vspace{6pt}
		
Next, we need a simple lemma of topological nature.
	\begin{Lem}\label{Lem:pi}
		Let $G$ be a Lie group acting on a smooth manifold $W$ such that $p:W\longrightarrow W/G$ is a principal $G$-bundle. Let $X$ be an open subset of $W$ and set $Z:=G\cdot X$. Then for every compact subset $K\subset p(Z)$ there exists a compact subset $\widetilde K\subset X$ such that $p^{-1}(K)=G\cdot \widetilde K$. 
	\end{Lem}
\Proof 
Let $K\subset p(Z)$ be a compact subset and $k\in K$. Since $p$ is a principal bundle, there exists an open neighborhood $U(k)$ of $k$ and a continuous section $s_k:U(k)\longrightarrow X$ with image in $X$. Then $K$ lies in a finite union of such neighborhoods. In particular, there are finitely many points $k_1, \dots, k_l \in K$ and relatively compact open subsets $A_j$ whose closure is contained in $U(k_j)$ such that $K=\bigcup\limits_{j=1}^l C_j$, where $C_j:=K\cap\overline{A_j}$. We set 
	$$
	\widetilde K:=\bigcup\limits_{j=1}^l s_{k_j}(C_j)\subset X.
	$$
By the continuity of $s_{k_j}$ the image of every $C_j$ is compact, therefore $\widetilde K$ is compact as a finite union of compact sets. As $p(\widetilde K)=K$ by construction, it follows that $p^{-1}(K)=G\cdot\widetilde K$, as wished.
\qed\vspace{6pt}
	
An important ingredient for the main result of this section is the following.
	\begin{Def}
		Let $f:T_\Omega\longrightarrow\RR$ be a smooth $H$-invariant strictly plurisubharmonic function and let $\pi:Z\longrightarrow Z/H^\CC$ be the quotient map. We say that $f$ is an exhaustion mod $H$ along the restriction $\pi|_{T_\Omega}$ of $\pi$ to $T_\Omega$ if for all compact subsets $K\subset Z/H^\CC$ and all $r\in\RR$, the set 
		$$
		\left(\pi^{-1}(K)\cap f^{-1}(-\infty,r]\right)/H\subset T_\Omega/H 
		$$ 
		is compact in the orbit space $T_\Omega/H$.
	\end{Def}
	Let $\rho:T_\Omega\longrightarrow\RR$ be the Kähler potential with respect to the Kähler form $\omega_\mathcal{B}$, introduced in Section \ref{twofour}.
	\begin{Prop}\label{Prop:exhaustion} 
		The Kähler potential $\rho$ is an exhaustion mod $H$ along $\pi|_{T_\Omega}$.
	\end{Prop}		
\Proof 
Let $K\subset Z/H^\CC$ be a compact subset and $r\in\RR$. As the quotient map $p:V^\CC\longrightarrow V^\CC/H^\CC$ is a principal bundle, we can apply Lemma \ref{Lem:pi}. This implies that $\pi^{-1}(K)=H^\CC+\widetilde K$ for a compact set $\widetilde K\subset T_\Omega$. Further, there exist compact sets $K_V\subset V$ and $K_\Omega\subset\Omega$ such that $\widetilde K\subset K_V+iK_\Omega$. Without loss of generality we can set $\widetilde K=K_V+iK_\Omega$. Then one has 
	$$
	\pi^{-1}(K)=H^\CC+(K_V+iK_\Omega)=H+K_V+i(H+K_\Omega).
	$$ 
Next, we observe that 
	$$
	\rho^{-1}(-\infty,r]=V+i\log(\varphi)^{-1}(-\infty,r].
	$$
The function $\log(\varphi)$ is strictly convex and goes to infinity on every sequence converging to the boundary of $\Omega$. Therefore, $\log(\varphi)^{-1}(-\infty,r]$ is closed in $V$, hence the set $(H+K_\Omega)\cap\log(\varphi)^{-1}(-\infty,r]$ is closed in $V$. In fact it is compact, as it is contained in the set $(H+K_\Omega)\cap\overline\Omega$ which is compact by Theorem \ref{Theo:compact}. This implies that the set 
	$$
	Q=K_V+i\left((H+K_\Omega)\cap\log(\varphi)^{-1}(-\infty,r]\right)
	$$
is compact. Finally, one has
	\begin{align*}
		\pi^{-1}(K)\cap\rho^{-1}(-\infty,r]=H+Q,
	\end{align*}
which implies the statement.
\qed
	\begin{Coro}\label{Coro:surj}
		For all $x\in T_\Omega$, we have $H^\CC\cdot x\cap\MM_H\neq\emptyset$. Hence, $T_\Omega\subset H^\CC\cdot\MM_H$, i.e., $Z=H^\CC\cdot\MM_H$.
	\end{Coro}
\Proof 
Recall from \cite{H}, p. 6, that
	$$
	\MM_H=\left\{x\in T_\Omega\ |\ x\text{ is a critical point of }\rho|_{H^\CC\cdot x\cap T_\Omega}\right\}.
	$$
Let $x=v+i\omega\in T_\Omega$. Proposition \ref{Prop:exhaustion} implies that $\rho$ has a minimum on 
    $$
    H^\CC\cdot x\cap T_\Omega=V+i((H+\omega)\cap\Omega).
    $$ 
In particular, this minimum is a critical point of the restriction of $\rho$ to $V+i((H+\omega)\cap\Omega)$, hence $x\in\MM_H$.
\qed
	\begin{Theo}\label{Theo:inclusiontrans}
		Let $T_\Omega$ be a homogeneous Siegel domain of the first kind and $H<V$ a connected subgroup of the translation group. Set $Z:=H^\CC+T_\Omega$. Then the map $\overline\iota:\MM_H/H\longrightarrow Z/H^\CC$ induced by the inclusion $\iota:\MM_{H}\longrightarrow Z$ is biholomorphic. In particular, the symplectic reduction $\MM_H/H$ is a Stein manifold biholomorphic to a Siegel domain of the first kind.
	\end{Theo}
\Proof 
Since $Z$ is open in $V^\CC$, the statement follows from Corollary \ref{Coro:inject} and Corollary \ref{Coro:surj}.
\qed 
	\begin{Rk}\textit{(Extension of Thm. \ref{Theo:inclusiontrans}.)}
		In Sect. 2.2 of \cite{Ku} we extend Theorem \ref{Theo:inclusiontrans} to connected subgroups of the group of real translations of Siegel domains of the second kind. For this we use a holomorphic $B$-equivariant embedding of $X$ into the Siegel upper half plane $D_S$ such that every connected subgroup $H$ of the group of real translations acts on $D_S$ as a connected subgroup of the group of real translations of $D_S$. Such an embedding is provided by \cite{Isym}. 
	\end{Rk}

\subsection{Symplectic reduction for subgroups which satisfy the Lie condition}\label{threethree}
Here we detect a class of subgroups of $B$ for which the symplectic reduction of $X$ turns out to be, as in the case of a translation subgroup, a Stein manifold. Let $H$ be a closed subgroup of $B$ and let $\mathcal{M}_H$ be the zero level set of $\mu_H$. As $H$ acts freely, $\MM_H$ is smooth and the momentum map $\mu_H$ is a submersion. For $x\in \MM_H$, consider the complex subspace $W_{x}:=T_{x}\MM_H\cap JT_{x}\MM_H$. One has for all $x\in \MM_H$,
	\begin{equation}\tag{$*$}\label{eq:star}
		\textup{Ker}(d\mu_H(x))=T_{x}\MM_H=(\mathfrak{h}\cdot x)^{\perp_{\omega_\mathcal{B}}}=\mathfrak{h}\cdot x\oplus W_{x}
    \end{equation}
and
	$$
	\dim_{\RR} \textup{Ker}(d\mu_H(x))=\dim_{\RR} X - \dim_{\RR} H = \dim_{\RR} W_x + \dim_{\RR} H.
	$$
	\begin{Def}\label{LieCondition}
		We say that the group $H$ satisfies the Lie condition if
		\begin{enumerate}[i)]
			\item $\mathcal{M}_H$ is connected,
			\item there exists $x_0$ in $\MM_H$ and a connected subgroup $S=\exp(\mathfrak{s})$ of $B$ (cf. Rk. \ref{Rk:connected}) such that $S\cdot x_0\subset\MM_H$ and $\mathfrak{s}\cdot x_0=W_{x_0}$.
		\end{enumerate}
	\end{Def}
	\begin{Lem}\label{Lem:HScap}
		Let $H$ be a closed subgroup of $B$ which satisfies the Lie condition. Then for all $x\in S\cdot x_0$, we have
		$$
		H\cdot x\cap S\cdot x=\left\{x\right\}.
		$$
	\end{Lem}
\Proof
We have $H\cap S=\{e\}$. Let $y\in H\cdot x\cap S\cdot x$. This implies $y=h\cdot x=s\cdot x$ for some $h\in H$ and $s\in S$ and therefore $s^{-1}h\cdot x=x$. Since $B$ acts freely, we have $h=s\in H\cap S=\left\{e\right\}$. This shows $y=x$. 
\qed\vspace{6pt}
		
The following proposition shows that if $H$ satisfies the Lie condition, then the $S$-orbit through $x_0$ is a closed complex submanifold of $X$. In particular, it is Stein. In fact, it is a homogeneous Siegel domain (see \cite{P}, pp.~52, 66). In the main theorem of this section we show that the symplectic reduction of $X$ is biholomorphic to $S\cdot x_0$.
	\begin{Prop}\label{Prop:complex}
		Let $X\subset\CC^d$ be a homogeneous Siegel domain and $H$ a closed subgroup of $B$ which satisfies the Lie condition.
		\begin{enumerate}[i)]
			\item The $S$-orbit through $x_0$ is a complex submanifold of $X$.
			\item The map 
			$$
			\Phi:H\times S\longrightarrow X,\ \Phi(h,s)=hs\cdot x_0
			$$ 
			is a diffeomorphism onto its image.
		\end{enumerate}
	\end{Prop}
\Proof 
i) The group $S$ acts on $X$ by holomorphic transformations and $S\cdot x_0=s\cdot (S\cdot x_0)$. This implies that the tangent space at each point $x\in S\cdot x_0$ is $J$-invariant. In particular, $S\cdot x_0$ is, as a submanifold of $X$, a complex manifold (see \cite{B}, p. 59).\vspace{6pt}
		
ii) The map $\Phi$ is a local diffeomorphism. Assume that $(h_1s_1)\cdot x_0=(h_2s_2)\cdot x_0$ for $h_1,h_2\in H$ and $s_1,s_2\in S$. Set $x:=s_1\cdot x_0$. Then one has
	$$
	(h_2^{-1}h_1)\cdot x=(s_2s_1^{-1})\cdot x.
	$$
Lemma \ref{Lem:HScap} implies that $h_1=h_2$ and $s_1=s_2$. Hence, $\Phi$ is injective. Indeed, note that for every $h\in H$ and $s\in S$ the real $H$-orbit through $hs\cdot x_0$ is transversal to the complex manifold $hS\cdot x_0$, whose tangent space at $hs\cdot x_0$ coincides with $W_{hs\cdot x_0}$ (cf. (\ref{eq:star})).
\qed\vspace{6pt}
		
The main result of this section follows from the following proposition.
	\begin{Prop}\label{Prop:Levi}
		Let $X\subset\CC^d$ be a homogeneous Siegel domain and $H$ a closed subgroup of $B$ which satisfies the Lie condition. Then the zero level set $\mathcal{M}_H$ is real analytic and Levi flat, i.e., $W:=T\mathcal{M}_H\cap JT\mathcal{M}_H$ is involutive and every leaf of the induced foliation is a complex manifold.
	\end{Prop}
\Proof 
The group $H$ acts real analytically and the Bergman metric is real analytic (see \cite{K}, p. 268). Therefore, $d\mu_H$ is real analytic. It follows that $\mu_H$ is real analytic and $\MM_H$ is a real analytic manifold. Let $W:=T\MM_H\cap JT\MM_H$. As $\MM_H$ is real analytic, $W$ is real analytic. Proposition~\ref{Prop:complex} implies that on the open set $HS\cdot x_0$ we have 
	$$
	W|_{HS\cdot x_0}=\bigsqcup\limits_{h\in H} T(hS\cdot x_0).
	$$
This implies that $W|_{HS\cdot x_0}$ is integrable and the set of all sections of $W|_{HS\cdot x_0}$ is a Lie algebra. In order to show that $W$ is involutive, we define the set
	$$
	M:=\left\{p\in\MM_H\ \left|
	\begin{array}{c}\exists\ \text{ open neighborhood }U \text{ of } p\text{ and}\\\text{ local sections }\{X_j\}\text{ of } W\\\text{ s.t. } W_x=\text{span}\{X_j|_{x}\}\text{ for all }x\in U\\\text{ and }[X_j,X_k]\text{ is a local section of $W$ for all } j\neq k\end{array}
	\right.\right\}.
	$$
The set $M$ is open by definition and contains $HS\cdot x_0$. As $\MM_H$ is connected by the definition of $H$, the statement will follow by showing that $M$ is closed. Let $\partial M$ be the topological boundary of $M$ and $p\in\partial M$. We claim that there exists an open neighborhood $U$ of $p$ and real analytic local sections $\{X_j\}$ of $W$ s.t. $W_x=\text{span}\{X_j|_{x}\}$ for all $x\in U$. On the open set $U\cap M$ the Lie bracket $[X_j,X_k]$ is a local section of $W$ for all $j\neq k$. As the local sections $X_j$ are real analytic, so are their Lie brackets. Thus, $[X_j,X_k]$ extends to all of $U$, implying $p\in M$.\\
It remains to prove the claim. For this observe that the bundle map $\pi_W:W\longrightarrow\MM_H$ is a real analytic submersion. By the constant rank theorem, there exist real analytic coordinates, on $W$ and on $\MM_H$ near $p$, such that $\pi_W$ reads as a linear projection with respect to these coordinates. In particular, there exist local real analytic sections of $\pi_W$ which proves the claim.
\qed

	\begin{Theo}\label{Theo:Quotient}
		Let $X\subset\CC^d$ be a homogeneous Siegel domain and $H$ a closed subgroup of $B$ which satisfies the Lie condition. Then
		\begin{enumerate}[i)]
			\item $\MM_H=HS\cdot x_0$,
			\item the restriction of the quotient map $\pi:\MM_{H}\longrightarrow\MM_H/H$ to $S\cdot x_0$ is a biholomorphic Kähler isometry.
		\end{enumerate}
		In particular, the symplectic reduction $\MM_H/H$ is Stein.
    \end{Theo}
\Proof 
i) As the set $HS\cdot x_0$ is open in $\MM_H$ by Proposition \ref{Prop:complex}, its projection $\pi(S\cdot x_0)$ is open in $\MM_H/H$. Assume by contradiction that there exists an element $y\in\MM_H$ such that $\pi(y)$ belongs to the boundary of $\pi(S\cdot x_0)$. Let $L_y$ be a local complex leaf of $W$ through $y$ (cf. Prop. \ref{Prop:Levi}). Since $T_y\MM_H=\mathfrak{h}\cdot y\oplus W_{y}$ (cf. (\ref{eq:star})), the map $H\times L_y\longrightarrow H\cdot L_y$, defined by $(h,x)\mapsto h\cdot x$, is a local diffeomorphism. It follows that the projection $\pi(L_y)$ is open. Hence, it intersects $\pi(S\cdot x_0)$. As a consequence, there exists $x\in L_y$ and $h\in H$ such that $x\in hS\cdot x_0$. Since $hS\cdot x_0$ is a closed (maximal) leaf, it follows that $L_y\subset hS\cdot x_0\subset HS\cdot x_0$. In particular $\pi(y)\in\pi(S\cdot x_0)$, which is a contradiction.\vspace{6pt}
		
ii) Statement (i) and Proposition \ref{Prop:complex} imply that the map $\Phi:H\times S\longrightarrow\MM_H$ is a diffeomorphism. Moreover, for all $x=s\cdot x_0\in S\cdot x_0$, we have 
	$$
	T_x(S\cdot x_0)=s_*(x_0)(W_{x_0})=W_x.
	$$
This implies that the restriction of $\pi$ to $S\cdot x_0$ is a Kähler isometry. 
\qed 
	\begin{Rk}(\textit{Testing the Lie condition.})
		Let $X$ be a homogeneous Siegel domain and $\mathfrak{b}$ the Lie algebra of $B=A\ltimes N$. Since $B$ acts simply transitively on $X$, we can identify its Lie algebra with the tangent space $T_{x_0}X$. In particular, this induces an integrable complex structure $J$ on $\mathfrak{b}$. Consider a root space decomposition of $\mathfrak{b}$ with respect to $\mathfrak{a}$. Let  $\mathfrak{h}<\mathfrak{b}$ be a Lie subalgebra defined as the direct sum of some root spaces. The structure theory of these root spaces together with an explicit description of $J$ (see \cite{P}, Ch. 2) allows us to identify candidate Lie subgroups $H=\exp(\mathfrak{h})$ for which a $J$-invariant Lie algebra and the associated group $S$ exist. Determining whether these candidates satisfy the Lie condition then reduces to checking that $S\cdot x_0\subset\MM_H$.
	\end{Rk}
	\begin{Ex}
		Let $T_\Omega=B\cdot x_0\subset V^\CC$ be a tube domain over a symmetric cone $\Omega\subset V$ and let $N_\Omega<G_\Omega$ be a maximal unipotent subgroup. Denote by $\mathfrak{b}=\mathfrak{a}\oplus\mathfrak{n}$ the Lie algebra of $B$ and set $r:=\dim_\RR\mathfrak{a}$. Choose linearly independent roots $\alpha_1,\dots,\alpha_r$ such that $J\mathfrak{a}$ is the direct sum of the corresponding root spaces, where $J$ denotes the complex structure on the tangent space $T_{x_0}T_\Omega\cong\mathfrak{b}$. Set $N_p:=\exp(J\mathfrak{a})$. The structure of the Lie algebra $\mathfrak{b}$, together with the computation of $\MM_{N_\Omega}$, implies that the group $N_\Omega$ satisfies the Lie condition with $S=AN_p$. In particular, $S\cdot x_0$ is biholomorphic to $(\Delta_1)^r\subset\CC^r$, where $\Delta_1$ denotes the unit disc.
	\end{Ex}
	\begin{Rk}(\textit{Action of the complexified group.})
		Let $X\subset\CC^d$ be a homogeneous Siegel domain and $H$ a closed subgroup of $B$ which satisfies the Lie condition. We have seen that the map $\Phi:H\times S\longrightarrow\MM_H$ is a diffeomorphism. It is natural to ask whether the map 
		$$
		H^\CC\times S\longrightarrow H^\CC S\cdot x_0
		$$ 
		is a biholomorphism and whether $H^\CC\times S$ is biholomorphic to $H^\CC\cdot X$, implying the latter manifold to be Stein. In the case where $X$ is the Lorentz tube (see Ex. \ref{Ex:lorentztube}) a positive answer to both questions is given in Sect.~3 of \cite{Ku} for closed subgroups of $N_\Omega$. Moreover, we develop methods that may enable a generalisation to the $m$-fold product $X^m$. This could constitute a first step towards a generalisation of the extended future tube conjecture proved in \cite{HSc}.
	\end{Rk}
	
\subsection*{Acknowledgements}
The author is very grateful to his supervisor Prof. Dr. Peter Heinzner for introducing him to this field of research and for numerous useful discussions and remarks. He is also grateful to the referees for many helpful comments and suggestions.

\bibliographystyle{amsalpha}
\bibliography{Hamiltonian_Actions_on_Homogeneous_Bounded_Domains}

@book{B,
	author = {Boggess, A.},
	title = {{CR} Manifolds and the Tangential {C}auchy-{R}iemann Complex},
	series = {Studies in Advanced Mathematics},
	publisher = {CRC Press},
	address = {Boca Raton},
	year = {1991},
}

@book{C,
	author = {Cartan, H.},
	title = {Sur les groupes de transformations analytiques},
	series = {(Exposés mathématiques IX.) volume 198 of Actual. scient. et industr.},
	publisher = {Hermann},
	address = {Paris},
	year = {1935},
}

@book{F,
	author = {Faraut, J. and Kor{\'a}nyi, A.},
	title = {Analysis on Symmetric Cones},
	series = {Oxford Mathematical Monographs},
	publisher = {Clarendon Press},
	year = {1994},
	address = {New York},
}

@article{GS,
	author = {Guillemin, V. and Sternberg, S.},
	title = {Geometric Quantization and Multiplicities of Group Representations},
	journal = {Invent. Math.},
	volume = {67},
	year = {1982},
	pages = {515--538},
}

@article{H,
	author = {Heinzner, P.},
	title = {The minimum principle from a {H}amiltonian point of view},
	journal = {Documenta Mathematica},
	year = {1998},
	volume = {3},
	pages = {1--14},
}

@article{HHL,
	author = {Heinzner, P. and Huckleberry, A. and Loose, F.},
	title = {K\"ahlerian extensions of the symplectic reduction},
	journal = {J. reine angew. Math.},
	year = {1994},
	volume = {455},
	pages = {123--140},
}

@article{HG,
	author = {Heinzner, P. and Schwarz, G. W.},
	title = {Cartan decomposition of the moment map},
	journal = {Math. Ann.},
	year = {2007},
	volume = {337},
	pages = {197--232},
}

@article{HK,
	author = {Heinzner, P. and Kutzschebauch, F.},
	title = {An equivariant version of {G}rauert's {O}ka principle},
	journal = {Invent. math.},
	year = {1995},
	volume = {119},
	pages = {317--346},
}

@article{HL,
	author = {Heinzner, P. and Loose, F.},
	title = {Reduction of complex {H}amiltonian {G}-spaces},
	journal = {Geometric and Functional Analysis},
	year = {1994},
	volume = {4},
	pages = {233--248},
}

@article{HSc,
	author = {Heinzner, P. and Schützdeller, P.},
	title = {The Extended Future Tube Conjecture for \textup{SO}(1,n)},
	journal = {Ann. Sc. Norm. Super. Pisa Cl. Sci. (5)},
	year = {2004},
	number = {1},
	volume = {III},
	pages = {39--52},
}

@article{HS,
	author = {Heinzner, P. and Stratmann, B.},
	title = {Invariant {K}\"ahler potentials and symplectic reduction},
	journal = {Ann. Sc. Norm. Super. Pisa Cl. Sci. (5)},
	year = {2023},
	number = {4},
	volume = {XXIV},
	pages = {2351--2402},
}

@article{Isym,
	author = "Ishi,H.",
	title = "\textit{On symplectic representations of normal j-algebras and their application to Xu's realizations of Siegel domains}",
	journal = "\textup{Differential Geom. Appl.}",
	year = "2006", 
	volume = "24",
	pages = "588-612",
}

@article{Kan,
	author = {Kaneyuki, S.},
	title = {On the automorphism groups of homogeneous bounded domains},
	journal = {J. Fac. Sci. Univ. Tokyo},
	volume = {14},
	year = {1967},
	pages = {89--130},
}

@book{Kn,
	author = {Knapp, A. W.},
	title = {Lie Groups Beyond an Introduction},
	series = {Progress in Mathematics},
	volume = {140},
	publisher = {Birkh\"auser},
	year = {2002},
	address = {Boston},
	edition = {2},
}

@article{Ku,
	author = {Kukol, M.},
	title = {Symplectic reduction for actions of selective groups on homogeneous bounded domains},
	journal = {PhD thesis, Ruhr-Universit\"at Bochum},
	year = {2025},
}

@article{Kur,
	author = {Kurtdere, A.},
	title = {Equivariant {K}\"ahlerian extensions of contact manifolds},
	journal = {J. reine angew. Math.},
	volume = {673},
	year = {2012},
	pages = {33--53},
}

@article{K,
	author = {Kobayashi, S.},
	title = {Geometry of Bounded Domains},
	journal = {Trans. Amer. Math. Soc.},
	volume = {92},
	number = {2},
	year = {1959},
	pages = {267--290},
}

@article{MW,
	author = {Marsden, J. and Weinstein, A.},
	title = {Reduction of symplectic manifolds with symmetry},
	journal = {Rep. Mathematical Phys.},
	year = {1974},
	volume = {5},
	number = {1},
	pages = {121--130},
}

@book{P,
	author = {Pyatetskii-Shapiro, I. I.},
	title = {Automorphic Functions and the Geometry of Classical Domains},
	publisher = {Gordon and Breach},
	year = {1969},
	address = {New York},
}

@book{R,
	author = {Rockafellar, R. T.},
	title = {Convex Analysis},
	series = {Princeton Mathematical Series},
	publisher = {Princeton University Press},
	year = {1970},
	address = {Princeton},
}

@article{Vin,
	author = {Vinberg, E. B.},
	title = {The {M}orozov-{B}orel theorem for real {L}ie groups},
	journal = {Dokl. Akad. Nauk SSSR},
	volume = {141},
	year = {1961},
	pages = {270--273},
}

@book{VinLie,
	editor = {Vinberg, E. B.},
	title = {{L}ie Groups and {L}ie Algebras {III}},
	series = {Encyclopaedia of Mathematical Sciences},
	volume = {41},
	year = {1994},
	publisher = {Springer},
	address = {Berlin},
}

\end{document}